\documentclass[letterpaper,10pt,oneside,twocolumn]{article}
\pdfoutput=1
\usepackage{graphicx}
\usepackage{geometry}
\usepackage[center]{titlesec}
\titleformat{\section}[hang]{\normalfont\normalsize\bfseries}{\thesection}{12pt}{\centering}%
\titleformat{\subsection}[display]{\normalfont\normalsize}{\thesubsection}{12pt}{\underline}%
\titleformat{\subsubsection}[runin]{\normalfont\normalsize}{\thesubsubsection}{12pt}{\underline}%
\newcommand{\PaperTitle}[1]{%
   \begin{center}%
      \begin{large}%
         \textbf {#1} \\%
      \end{large}%
   \end{center}%
}%
\newcommand{\AuthorList}[1]{%
   \begin{center}%
      {#1} \\%
   \end{center}%
}%
\newcommand{\AuthorAffiliations}[1]{%
   \begin{center}%
      {#1}%
   \end{center}%
}%
\newcommand{\Keywords}[1]{%
   \begin{center}%
      Keywords: {#1}%
   \end{center}%
}%
\begin{document}%
\twocolumn[%
%
%
\PaperTitle{An analysis of non-classical austenite-martensite interfaces in CuAlNi}%
%
%
\AuthorList{J. M. Ball$^1$, K. Koumatos$^1$, H. Seiner$^{2,3}$}%
\AuthorAffiliations{$^1$Mathematical Institute, University of Oxford, 24-29 St. Giles', Oxford OX1 3LB\\ $^2$Institute of Thermomechanics ASCR, Dolej\v skova 5, 182 00 Prague 8, Czech Republic\\ $^3$Faculty of Nuclear Sciences and Physical Engineering, CTU in Prague, Trojanova 13, 120 00 Prague 2, Czech Republic}%
\Keywords{martensitic transformations, non-classical austenite-martensite interfaces, cubic-to-orthorhombic}%
]%
%
%
\section{Abstract}  Ball \& Carstensen \cite{1,2}, theoretically investigated the possibility of the occurrence of non-classical austenite-martensite interfaces and studied the cubic-to-tetragonal case extensively. Here, we aim to present an analysis of such interfaces recently observed by Seiner et al. \cite{3} in CuAlNi single crystals, undergoing a cubic-to-orthorhombic transition. We show that they can be described by the non-linear elasticity model for martensitic transformations and we make some predictions regarding the volume fractions of the martensitic variants involved, as well as the habit plane normals.%
\section{Introduction} A classical austenite-martensite interface is one in which the undistorted austenite meets a simple laminate of martensite and these have been broadly studied. In recent years, a theory of martensitic transformations has been developed which allows austenite-martensite interfaces to occur in which the microstructure of the martensite is more complicated; these are referred to as non-classical interfaces and have hitherto not been systematically observed. In this non-linear elasticity model, in which interfacial energy is neglected, microstructures are identified as limits of minimizing sequences for the total free energy
\begin{equation}
	I_{\theta}(y)=\int_{\Omega}\varphi(\nabla y(x),\theta)dx.
	\label{eq:1}
\end{equation}
Here, $\Omega$ represents the reference configuration of undistorted austenite at the transformation temperature $\theta_{c}$ and $y(x)$ denotes the deformed position of particle $x\in\Omega$. The free-energy function $\varphi(F,\theta)$ depends on the deformation gradient $F$ and the temperature $\theta$. By frame indifference, $\varphi(RF,\theta)=\varphi(F,\theta)$ for all $F$, $\theta$ and for all rotations $R$; that is for all $3\times3$ matrices in $SO(3)=\left\{R:R^{T}R=\mathbf{1}, \det{R}=1\right\}$. By adding an appropriate function of $\theta$, we may assume that $\min_{F} \varphi(F,\theta)=0$. At $\theta_{c}$, the energy wells of the free-energy function are given by $SO(3)$ for the undistorted austenite and $SO(3)U_{i}$, $i=1,...,N$, for the $N$ distinct variants of martensite, where each $U_{i}$ is a positive definite, symmetric matrix. Hence, $\varphi(F,\theta_{c})\geq0$ with equality if and only if $F\in SO(3)\cup\bigcup^{N}_{i=1}SO(3)U_{i}$. For $\theta<\theta_{c}$, the martensite wells, $\bigcup^{N}_{i=1}SO(3)U_{i}$, minimize $\varphi$, whereas for $\theta>\theta_{c}$ the minimum is given by the austenite well, $SO(3)\alpha(\theta)\mathbf{1}$, where $\alpha(\theta)$ is the thermal expansion coefficient of the austenite and $\alpha(\theta_{c})=1$.

A non-classical planar austenite-martensite interface $\left\{x\cdot m=k\right\}$ corresponds to a choice of habit-plane normal $m$ such that there exists a sequence of deformations $y^{(j)}$ for which $I_{\theta_{c}}(y^{(j)})\longrightarrow0$ as $j\longrightarrow\infty$. We require that, for $x\cdot m<k$, the values of the deformation gradient $\nabla y^{(j)}$ tend to $SO(3)$, i.e. $y^{(j)}$ corresponds to the austenite phase; without loss of generality, this is equivalent to $\nabla y^{(j)}\left(x\right) \longrightarrow\mathbf{1}$ except possibly for a set of zero volume. For $x\cdot m>k$, we require that, as $j\longrightarrow\infty$, $\nabla y^{(j)}$ tends in a suitable way to the set $K=\bigcup^{N}_{i=1} SO(3)U_{i}$ of martensite energy wells. In fact, we require that the Young measure $(\nu_{x})_{x\in\Omega}$ of $\nabla y^{(j)}$ is supported in the set $K$ (for details see e.g. \cite{6}).

From now on, we shall make the assumption that the martensitic microstructure is homogeneous; that is, for $x\cdot m>k$, the macroscopic deformation gradient $F=\nabla y(x)$ is independent of $x$. The set of all such matrices $F$ is called the \textit{quasiconvexification} of $K$ and is denoted by $K^{qc}$. It can be shown that $K^{qc}$ is also the set of $F$ such that there exists a sequence of deformations $z^{(j)}$ with $z^{(j)}(x)=Fx$ on the boundary of $\Omega$, $\partial\Omega$, and $\nabla z^{(j)}(x)\longrightarrow K$ in the above sense.

If we know $K^{qc}$ for a given set of martensite wells, to ensure geometric compatibility, we need to examine whether it is possible to establish a rank-one connection between $SO(3)$ and $K^{qc}$; that is, we need to find vectors $b$ and $m$ such that
\begin{equation}
	\mathbf{1}+b\otimes m\in K^{qc},
	\label{eq:2}
\end{equation}
where, by frame indifference, we have chosen the identity matrix $\mathbf{1}$ to represent the austenite energy well. Unfortunately, we only have a characterization of $K^{qc}$ for two martensitic wells ($N=2$), i.e. when $K=SO(3)U_{1}\cup SO(3)U_{2}$. In this case, any $F\in K^{qc}$ can be obtained as the macroscopic deformation gradient of a double laminate (see \cite{13}). Even though $K^{qc}$ is unknown in the case of three tetragonal wells, Ball \& Carstensen \cite{2} were able, using the two-well calculation, to characterize exactly the values of the deformation parameters for a cubic-to-tetragonal transformation which permit non-classical interfaces. They also presented results for a cubic-to-orthorhombic transformation, though not of the type occurring in CuAlNi.

Given a matrix $F$, the question as to whether we can solve the equation 
\begin{equation}
	\mathbf{1}+b\otimes m\in SO(3)F
	\label{eq:3}
\end{equation}
for vectors $b$ and $m$ is answered by the following lemma in \cite{8,12}.

\textbf{Lemma 1.} Let $F$ be a non-singular matrix that is not a rotation. Then the wells $SO(3)$ and $SO(3)F$ are rank-one connected if and only if the middle eigenvalue of $F^{T}F$ is 1. Then $\mathbf{1}+b\otimes m\in SO(3)F$ for some $b$ if and only if m is a non-vanishing multiple of one of the two vectors, $\sqrt{1-\lambda_{1}}e_{1}\pm\sqrt{\lambda_{3}-1}e_{3}$, where $0\leq\lambda_{1}\leq1\leq\lambda_{3}$ are the three eigenvalues of $F^{T}F$ with corresponding eigenvectors $e_{1}$, $e_{2}$, $e_{3}$.

Having outlined a brief description of the model, we proceed to the experimental observations on a CuAlNi single crystal.
\section{Experimental Observations} 
The first micrographs of interfaces between austenite and crossing twins were obtained by Seiner et al. \cite{4}, who documented that such interfaces can form during the shape recovery process of CuAlNi single crystals. However, these interfaces were only \emph{weakly non-classical}, i.e. they were classical interfaces weakly disturbed by a negligible volume fraction of compound twins intersecting the first-order laminate of the Type-II twins. These observations motivated the authors towards more intensive research in this field -- the experiment was improved in order to increase the volume fraction of the compound twins in the martensitic microstructure. The resulting experimental procedure is briefly outlined below, and will appear in more details in \cite{3}.

\begin{figure}[h]
	\centering
		\includegraphics[width=0.4\textwidth]{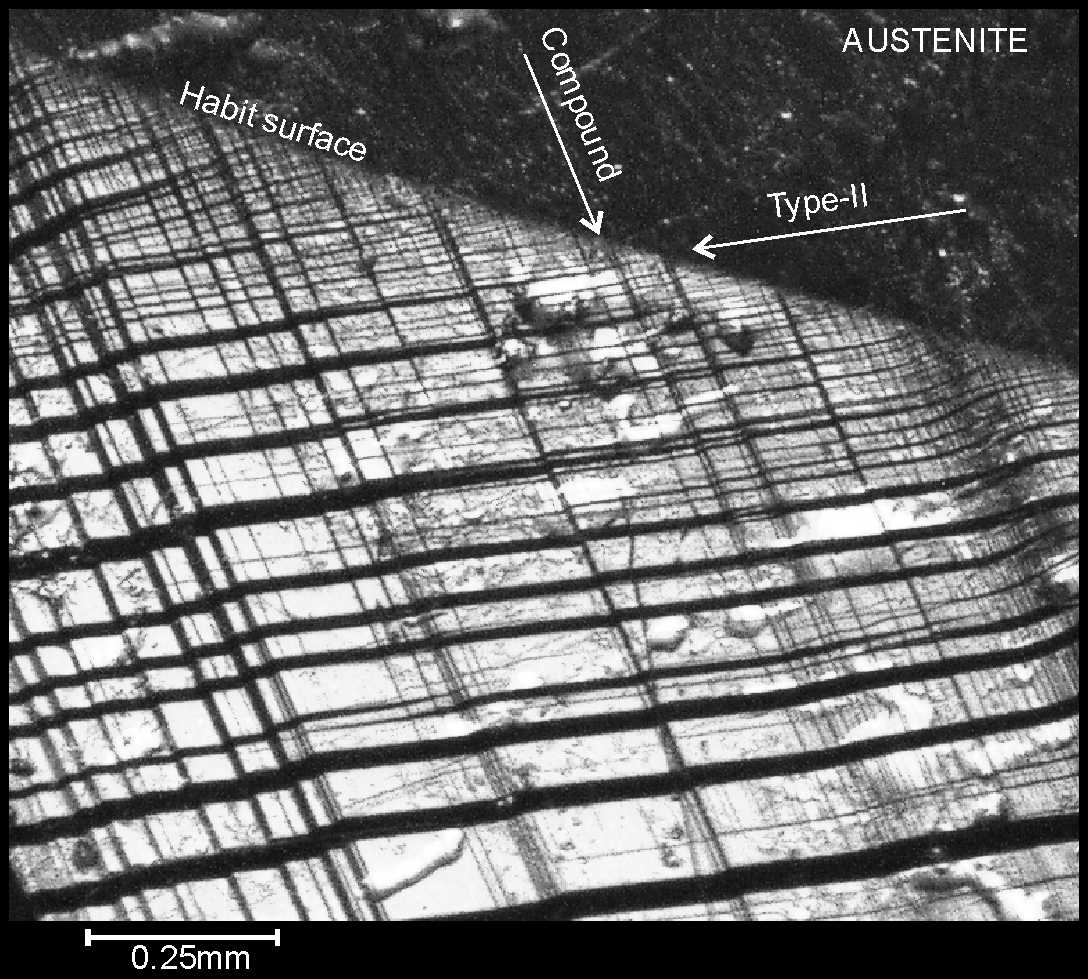}
	\caption{Optical micrograph of a non-classical interface between austenite and a martensitic microstructure. The arrows indicate the orientations of twinning planes of Type-II and compound twinning systems.}
	\label{fig:Non-calssical Interface}
\end{figure}

\begin{figure}[h]
	\centering
		\includegraphics[width=0.4\textwidth]{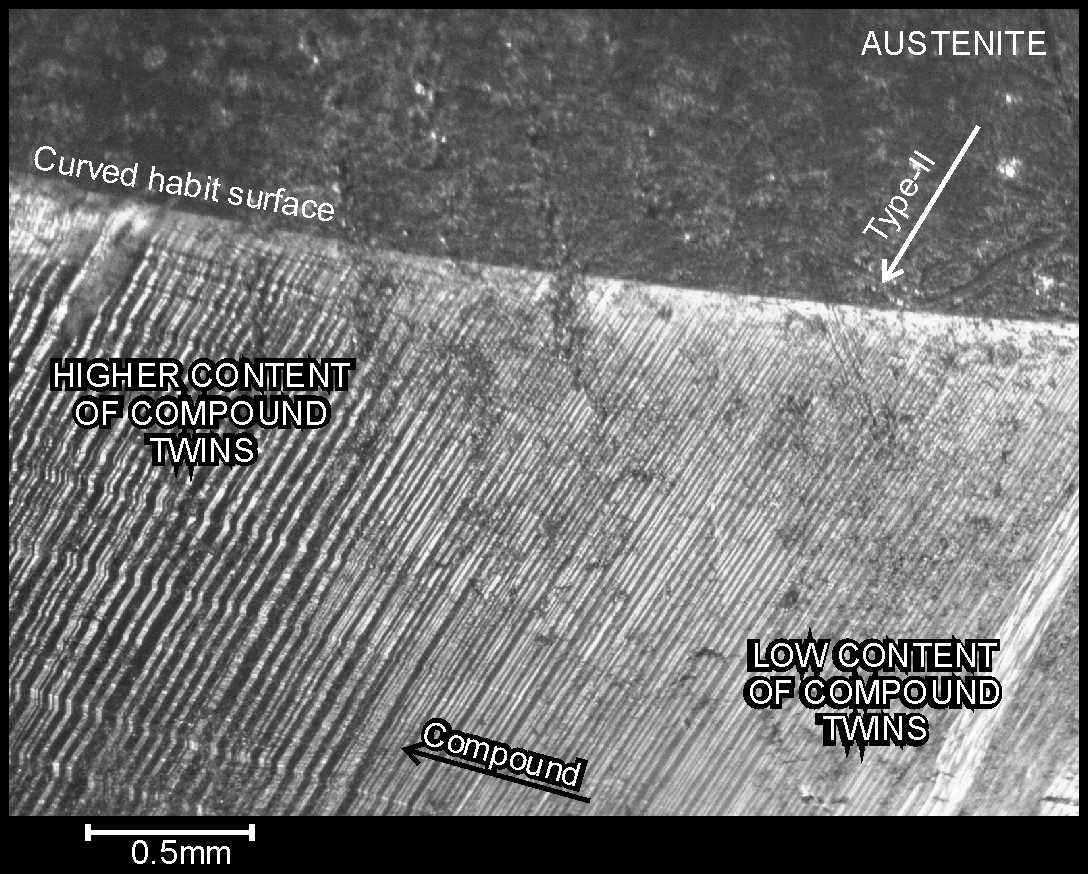}
	\caption{Curved interface between crossing twins and austenite resulting from the inhomogeneity of compound twinning. (Optical microscopy.)}
	\label{Curved Interface}
\end{figure}

The specimen examined in this case was a 3.9$\times$3.8$\times$4.2mm rectangular parallelepiped of the austenitic phase of CuAlNi, cut from a single crystal of this alloy such that the normals to the specimen faces had approximately the principal crystallographic directions $\langle$100$\rangle$. The original single crystal was grown by a Bridgman method at the Institute of Physics ASCR in Prague.
The specimen was subjected to the following sequence of mechanical and thermal loadings (see Fig.~\ref{nonc}):
\begin{figure}[ht]
\centering
\includegraphics[width=0.35\textwidth]{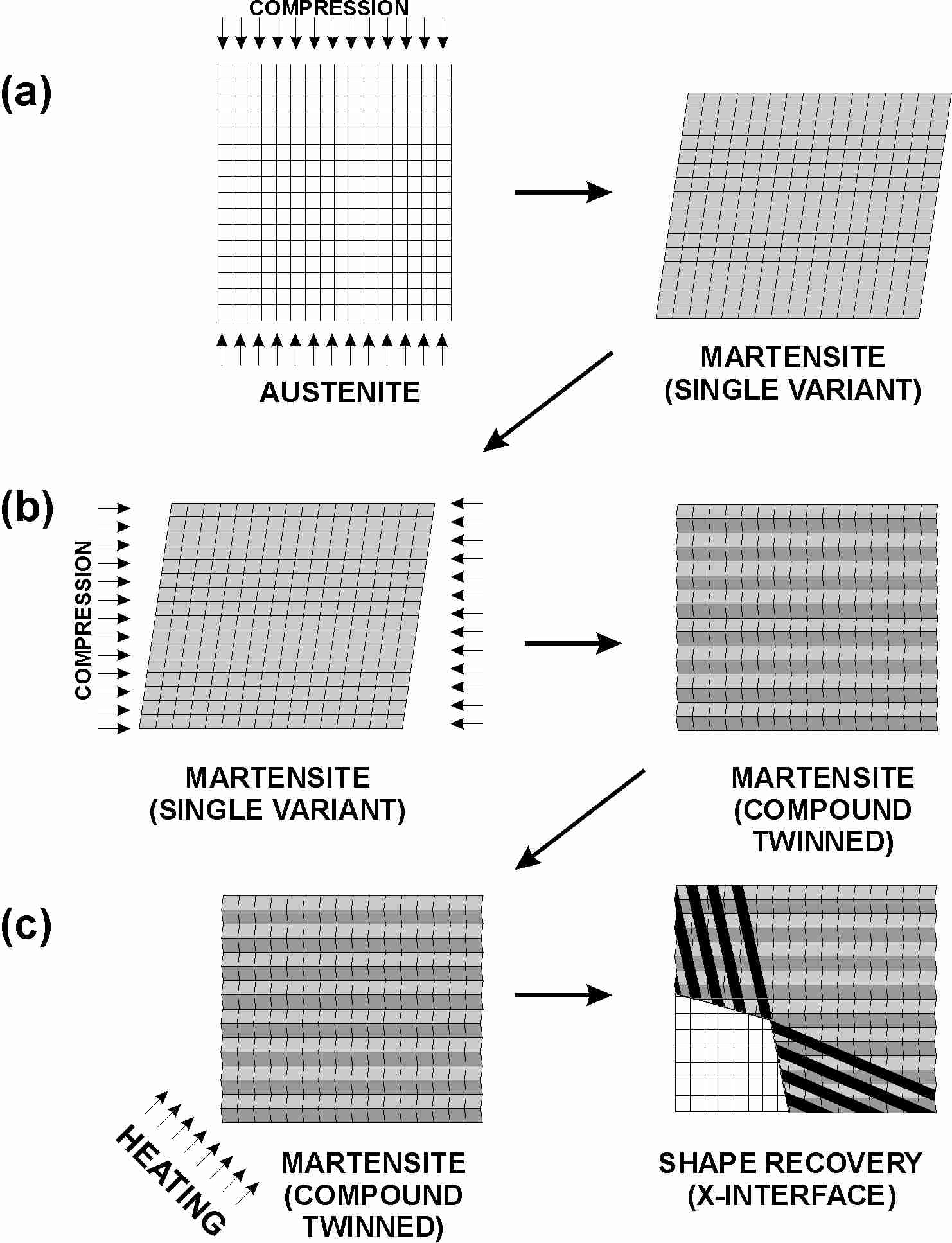}
\caption{Outline of the experimental procedure.}
\label{nonc}
\end{figure}
\begin{itemize}
\item[(a)]{At room temperature, the specimen was transformed into a single variant of 2H martensite by applying uniaxial compression (in a bench vice). Due to the effect called mechanical stabilization of martensite, the specimen did not return to austenite after unloading, but remained as a single variant of martensite.}
\item[(b)]{The specimen was rotated by 90$^\circ$ and uniaxial compression was applied again. In this case, the loading induced the reorientation of martensite into another variant via compound twinning (for an analysis of the relation between the direction of applied compression and the activated twinning systems in CuAlNi, see \cite{5}). The reorientation was not fully completed. Instead, the loading was interrupted at the moment when the specimen contained comparable volume fractions of both variants. By such a procedure, we obtained a finely compound twinned specimen.}
\item[(c)]{The finely compound twinned specimen was heated from one side using a gas lighter, which induced the shape recovery process, i.e. the thermally driven return of the specimen into austenite. As the compound twins cannot form any compatible interface with austenite, the transition was achieved by formation of an \emph{interfacial microstructure}, which ensured a compatible connection between the mechanically stabilized martensite (the compound twins) and austenite. This interfacial microstructure was formed by Type-II twins crossing the original compound microstructure and getting arranged into a so called X-interface (for more details of formation of X-interfaces in CuAlNi see \cite{4}, for the theoretical analysis of this microstructure, see \cite{7}). By removing the heating in the middle of the course of the transition, the interfacial microstructure was stopped, and the non-classical interfaces between austenite and the two mutually crossing systems of twins (compound and Type-II) were observable by optical microscopy. An example is given in Fig.~\ref{fig:Non-calssical Interface}.}
\end{itemize}
This procedure was repeated several times, which enabled the capturing of several optical micrographs of the non-classical interfaces. An interesting observation was that, as shown in Fig.~\ref{Curved Interface}, the interface between the austenite and the system of crossing twins was never exactly planar, but rather slightly curved. This results from the fact that the pattern of compound twins induced in the specimen by the uniaxial compression in stage (b) of the experimental procedure is never exactly homogeneous. With varying volume fraction of the compound twins, the admissible orientation of the habit plane varies as well, as will be shown in the theoretical analysis given in the following section. However, a more complete analysis of the curved interface remains to be done.
\section{Analysis of Non-Classical Interfaces} 
\begin{figure}[h]
	\centering
		\includegraphics[scale=0.4]{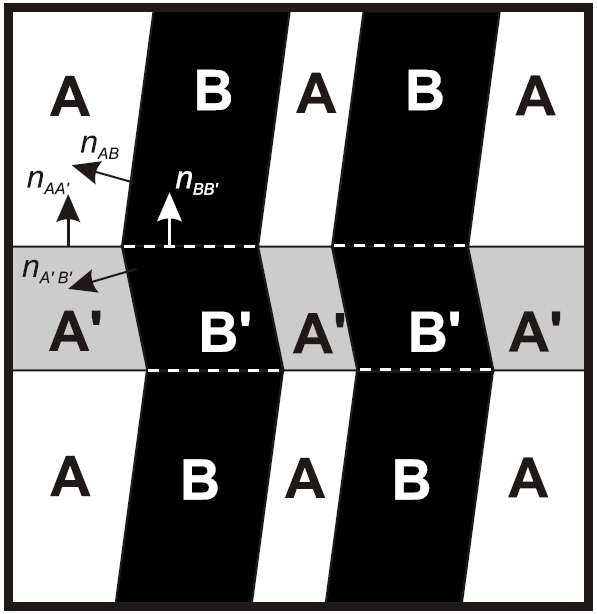}
	\caption{Parallelogram microstructure}
	\label{fig:para2}
\end{figure}

In this section, we present an analysis for the above non-classical austenite-martensite interface via the non-linear elasticity model for martensitic transformations. The martensitic region on one side of the observed interface consists of twin crossings involving four martensitic variants. Since the quasiconvexification of more than two wells is unknown, we are not able to analyze all the possibilities for non-classical interfaces in CuAlNi. However, we can do this in the twin crossing case. Firstly, we shall concentrate on the martensite phase and try to construct the microstructure consisting of compound and Type-II twin crossings, as shown in Fig.~\ref{fig:para2}.

Let $U_{A}$ and $U_{B}$ be two martensitic variants able to form a Type-II twin and $U_{A'}$, $U_{B'}$ the respective compound counterparts. Clearly $U_{A'}$ and $U_{B'}$ can also form a Type-II twin and we proceed by writing the compatibility equations for the parallelogram microstructure. These are:
\begin{equation}
R_{AB}U_{B}-U_{A}=b_{AB}\otimes n_{AB}
\label{eq:4}
\end{equation}
\begin{equation}
R_{A'B'}U_{B'}-U_{A'}=b_{A'B'}\otimes n_{A'B'}
\label{eq:5}
\end{equation}
\begin{equation}
R_{AA'}U_{A'}-U_{A}=b_{AA'}\otimes n_{AA'}
\label{eq:6}
\end{equation}
\begin{equation}
R_{BB'}U_{B'}-U_{B}=b_{BB'}\otimes n_{BB'}
\label{eq:7}
\end{equation}
where $b_{IJ}$ and $n_{IJ}$ are, respectively, the shearing vector and the normal to the twinning plane for the system of variants $U_{I}$ and $U_{J}$. Also, $R_{IJ}$ denotes the mutual rotation between variants $U_{I}$ and $U_{J}$.

If the above compatibility equations hold, a necessary and sufficient condition for the entire parallelogram microstructure to be compatible is that
\begin{equation}
R_{AB}R_{BB'}=R_{AA'}R_{A'B'}.
\label{eq:8}
\end{equation}
Equation (\ref{eq:8}) is necessary, since both sides describe the mutual rotation between $U_{A}$ and $U_{B'}$, and sufficiency follows by showing that (\ref{eq:4})-(\ref{eq:7}) imply that the normals $n_{AB}$, $n_{A'B'}$, $n_{AA'}$ and $n_{BB'}$ are coplanar.

Let
\begin{equation}
	M_{AB}=(1-\lambda)U_{A}+\lambda R_{AB}U_{B}
	\label{eq:9}
\end{equation}
\begin{equation}
	M_{A'B'}=(1-\lambda)U_{A'}+\lambda R_{A'B'}U_{B'}
	\label{eq:10}
\end{equation}
which represent the macroscopic deformation gradients corresponding to the Type-II structures. From \cite{11}, there are solutions of the twinning relations (\ref{eq:6}) and (\ref{eq:7}) with
\begin{equation}
n_{AA'}=n_{BB'}
\label{eq:11}
\end{equation}
and we make this choice in accordance with the experimental observations. Then the geometry of the parallelogram structure requires $\lambda$ to be the same in (\ref{eq:9}) and (\ref{eq:10}).

At this stage we note that $\det M_{A'B'}=\det M_{AB}$. To see this, consider equation (\ref{eq:4}). We get that
\[R_{AB}U_{B}=U_{A}(1+U^{-1}_{A}b_{AB}\otimes n_{AB}).\]
Therefore, by taking determinants, we obtain
\[\det U_{B}=\det U_{A}(1+U^{-1}_{A}b_{AB}\cdot n_{AB}).\]
Clearly $\det U_{B}=\det U_{A}$ and thus
\begin{equation}
	U^{-1}_{A}b_{AB}\cdot n_{AB}=0.
	\label{eq:12}
\end{equation}
Hence, for $M_{AB}$, we have
\[\det M_{AB}=\det U_{A}.\]
Similarly, $\det M_{A'B'}=\det U_{A'}=\det U_{A}$ and we deduce that $\det M_{A'B'}=\det M_{AB}$. This result has a simple physical interpretation. Since the determinants have the meaning of volume change and the microstructures $M_{AB}$ and $M_{A'B'}$ are just mixtures of variously rotated single variants, the volume must remain the same.

Now, assume that the volume fraction of variant $U_{A'}$ in $U_{A}$, as well as that of $U_{B'}$ in $U_{B}$, is $\Lambda$. Then, the macroscopic deformation gradient of the entire microstructure is
\begin{equation}
	M=(1-\Lambda)M_{AB}+\Lambda R_{AA'}M_{A'B'}.
	\label{eq:13}
\end{equation}
Whether the interface between the parallelogram microstructure and the austenite can be formed is a question of finding $\lambda$ and $\Lambda$ such that the middle eigenvalue of $M^{T}M$ is equal to 1, which is equivalent to finding a rank-one connection between the austenite and $SO(3)M$. In particular, by Lemma 1,
\begin{equation}
	\det(M^{T}M-\mathbf{1})=0.
	\label{eq:14}
\end{equation}
\subsection{Forming the Interface} By manipulating the compatibility equations and making use of equations (\ref{eq:8}) and (\ref{eq:11}) we can see that
\[(b_{AA'}-R_{AB}b_{BB'})\otimes n_{AA'}-b_{AB}\otimes n_{AB}+R_{AA'}b_{A'B'}\otimes n_{A'B'}\]\[=0.\]
All normals cannot be parallel to each other and hence, using a result in \cite{7}, we deduce that
\begin{equation}
	b_{AA'}-R_{AB}b_{BB'}\parallel b_{AB}\parallel R_{AA'}b_{A'B'}.
	\label{eq:15}
\end{equation}
Hence, there exists some $\eta$ such that
\begin{equation}
\eta b_{AB}=b_{AA'}-R_{AB}b_{BB'}.
\label{eq:16}
\end{equation}
By identifying the rotation $R_{AB}$ and using the formulae for the twinning solutions in Result 5.2 of \cite{10}, as well as the relations between them \cite{11}, we obtain that 
\[\eta=\frac{2}{\left|b_{AB}\right|^{2}}b_{AA'}\cdot b_{AB}.\]
Our goal is to deduce an expression for $\det(M^{T}M-1)$ which will enable us to find solutions of (\ref{eq:14}) for $\lambda, \Lambda\in(0,1)$. Towards this end, we shall seek a rank-one connection between $M_{AB}$ and $M_{A'B'}$. Indeed, using (\ref{eq:8}) and equations (\ref{eq:4})-(\ref{eq:7}),
\[R_{AA'}M_{A'B'}-M_{AB}=(1-\lambda)b_{AA'}\otimes n_{AA'}+\lambda R_{AB}b_{BB'}\otimes n_{BB'}\]
Having that $n_{AA'}=n_{BB'}$ and using (\ref{eq:16}) we get
\begin{equation}
	R_{AA'}M_{A'B'}-M_{AB}=b_{*}\otimes n_{AA'},
\label{eq:17}
\end{equation}
where $b_{*}=b_{AA'}-\lambda\eta b_{AB}$

Combining equation (\ref{eq:17}) with (\ref{eq:13}) we get that
\[M=M_{AB}+\Lambda b_{*}\otimes n_{AA'}\]
and using the expression for $M_{AB}$ we obtain
\begin{equation}
	M=U_{A}+\lambda b_{AB}\otimes n_{AB}+\Lambda b_{*}\otimes n_{AA'}.
	\label{eq:18}
\end{equation}
Then,
\begin{equation}
	\det(M^{T}M-\mathbf{1})=\det M^{T}\det(M-M^{-T}).
	\label{eq:19}
\end{equation}
We shall first calculate $\det M^{T}=\det M$. Hence we have,
\[\det M=\det M_{AB}\det(1+\Lambda M^{-1}_{AB}b_{*}\otimes n_{AA'}).\]
However, from (\ref{eq:17})
\[R_{AA'}M_{A'B'}=M_{AB}(1+M^{-1}_{AB}b_{*}\otimes n_{AA'})\]
and taking determinants, we see that
\[\det M_{A'B'}=\det M_{AB}(1+M^{-1}_{AB}b_{*}\cdot n_{AA'}).\]
Since $\det M_{A'B'}=\det M_{AB}$ we get that
\begin{equation}
	M^{-1}_{AB}b_{*}\cdot n_{AA'}=0.
	\label{eq:20}
\end{equation}
From (\ref{eq:19}), we have
\[\det(M^{T}M-\mathbf{1})=\det U_{A}\det(M-M^{-T})\]
and it remains to calculate $M-M^{-T}$.

Making use of the expression for $M$ and (\ref{eq:20})
\[M^{-T}=M^{-T}_{AB}-\Lambda M^{-T}_{AB}n_{AA'}\otimes M^{-1}_{AB}b_{*}\]
and from (\ref{eq:9}) and (\ref{eq:12})
\[M^{-1}_{AB}=U^{-1}_{A}-\lambda U^{-1}_{A}b_{AB}\otimes U^{-1}_{A}n_{AB}.\]
Combining the last two equations and after some calculations we obtain that
\begin{equation}
	M^{-T}=U^{-1}_{A}-\lambda U^{-1}_{A}n_{AB}\otimes U^{-1}_{A}b_{AB}
	\label{eq:21}
\end{equation}
\[-\Lambda U^{-1}_{A}n_{AA'}\otimes U^{-1}_{A}b_{AA'}\]
\[+\Lambda\lambda U^{-1}_{A}n_{AA'}\otimes U^{-1}_{A}b_{AB}(U^{-1}_{A}n_{AB}\cdot b_{AA'}+\eta).\]
Here we have already used the fact that
\begin{equation}
	U^{-1}_{A}n_{AA'}\cdot b_{AB}=0.
	\label{eq:22}
\end{equation}
The above result requires some effort and is provided by investigating the axes of the rotations in the austenitic point group along with Result 5.2 in \cite{10}.

Using the expressions for $M$ and $M^{-T}$ we get
\begin{equation}
	\det(M^{T}M-\mathbf{1})=\det(A_{0}+\lambda A_{1}+\Lambda A_{2}+\lambda\Lambda A_{3}),
	\label{eq:23}
\end{equation}
where
\begin{eqnarray}
\label{eq:24}
A_{0}&=&U^{2}_{A}-1\\
\label{eq:25}
A_{1}&=&U_{A}b_{AB}\otimes n_{AB}+n_{AB}\otimes U^{-1}_{A}b_{AB}\\
\label{eq:26}
A_{2}&=&U_{A}b_{AA'}\otimes n_{AA'}+n_{AA'}\otimes U^{-1}_{A}b_{AA'}\\
\label{eq:27}
A_{3}&=&-(U^{-1}_{A}n_{AB}\cdot b_{AA'}+\eta)n_{AA'}\otimes U^{-1}_{A}b_{AB}\\
         &&-\eta U_{A}b_{AB}\otimes n_{AA'}\nonumber .
\end{eqnarray}
It is trivial to observe that, for fixed $\lambda$, the expression multiplying $\Lambda$ is of rank 2 and similarly, for fixed $\Lambda$, the expression multiplying $\lambda$ is of the same rank. Hence, the determinant of $M-M^{-T}$ will only contain terms with $\lambda$ and $\Lambda$ in powers not greater than two. Letting $g(\lambda,\Lambda):=\det(M^{T}M-\mathbf{1})$ we deduce that
\begin{equation}
\label{eq:28}
g(\lambda,\Lambda)=\alpha+\beta\lambda+\gamma\Lambda+a\lambda^{2}+b\Lambda^{2}+c\lambda\Lambda+d\lambda\Lambda^{2}
\end{equation}
\[+e\lambda^{2}\Lambda+f\lambda^{2}\Lambda^{2}.\]
Using the relations between the twinning solutions for different martensitic variants \cite{11}, the above expression simplifies further by noting that for all $\Lambda\in(0,1)$
\begin{equation}
	g(0,\Lambda)=g(1,\Lambda)
	\label{eq:29}
\end{equation}
and, respectively, for all $\lambda\in(0,1)$
\begin{equation}
	g(\lambda,0)=g(\lambda,1).
	\label{eq:30}
\end{equation}
Now, from equations (\ref{eq:28})-(\ref{eq:30}) we obtain the following form for $g$;
\begin{equation}
	 g(\lambda,\Lambda)=a_{0}+a_{1}(\lambda^{2}-\lambda)+a_{2}(\Lambda^{2}-\Lambda)+a_{3}(\lambda^{2}-\lambda)(\Lambda^{2}-\Lambda).
	 \label{eq:31}
\end{equation}
From (\ref{eq:23}) and (\ref{eq:31}) it becomes easy to specify the coefficients $a_{0}, a_{1}, a_{2}, a_{3}$ in the expression for $g(\lambda,\Lambda)$. Firstly,
\begin{equation}
	a_{0}=g(0,0)=\det(U^{2}_{A}-1).
	\label{eq:32}
\end{equation}
Noting that $a_{1}=-\frac{\partial g}{\partial \lambda}(0,0)$ and $a_{2}=-\frac{\partial g}{\partial \Lambda}(0,0)$, we deduce, after a simple calculation, that
\begin{equation}
	a_{1}=-2b_{AB}\cdot U_{A}\textrm{cof}(U^{2}_{A}-1)n_{AB}
	\label{eq:33}
\end{equation}
and
\begin{equation}
	a_{2}=-2b_{AA'}\cdot U_{A}\textrm{cof}(U^{2}_{A}-1)n_{AA'}.
	\label{eq:34}
\end{equation}
As for the last coefficient, $a_{3}=4(\frac{\partial g}{\partial\lambda}(0,\frac{1}{2})+a_{1})$ and the expression gets more complicated reducing to
\begin{equation}
a_{3}=4\mathrm{cof}(A_{0}+\frac{1}{2}A_{2})\cdot(A_{1}+\frac{1}{2}A_{3})+4a_{1}.
\label{eq:35}
\end{equation}
For $\Lambda=0$, we only have the Type-II structure and we get
\[g(\lambda,\Lambda)=a_{0}+a_{1}(\lambda^{2}-\lambda).\]
This becomes zero for $\lambda^{2}-\lambda=-\frac{a_{0}}{a_{1}}$; call this $\lambda=\lambda^{*}\in\left(0,\frac{1}{2}\right] $, which agrees with the value obtained, for example, by Ball \& James in \cite{8}.

Setting $g(\lambda,\Lambda)=0$ and solving for $\lambda$ we deduce that, as long as
\begin{equation}
\frac{a_{1}}{a_{3}}\notin[0,\frac{1}{4}]
\label{eq:36}
\end{equation}
then
\begin{equation}
\lambda^2-\lambda=-\frac{a_{0}+a_{2}(\Lambda^2-\Lambda)}{a_{1}+a_{3}(\Lambda^2-\Lambda)}.
\label{eq:37}
\end{equation}
For $\Lambda=0$ and $\Lambda=1$ the above equation admits two solutions, namely $\lambda^{*}$ and $1-\lambda^{*}$, and we see that branches starting from $(\lambda^{*},0)$ and $(\lambda^{*},1)$ are created consisting of values of $(\lambda,\Lambda)$ that make the non-classical interface possible. Provided that 
\begin{equation}
a_{0}a_{3}\neq a_{1}a_{2}
\label{eq:38}
\end{equation}
a necessary and sufficient condition for these to meet, i.e. for equation (\ref{eq:37}) to have solutions for all $\Lambda\in\left[0,1\right] $, is that
\begin{equation}
0\leq\frac{4a_{0}-a_{2}}{4a_{1}-a_{3}}\leq\frac{1}{4}.
\label{eq:39}
\end{equation}
If, in addition,
\begin{equation}
\frac{4a_{0}-a_{2}}{4a_{1}-a_{3}}\neq\frac{1}{4}
\label{eq:40}
\end{equation}\\
then the solutions will be distinct.

Otherwise, if $\alpha_{0}\alpha_{3}=\alpha_{1}\alpha_{2}$, $g$ simplifies even further and we obtain that, for all $\Lambda\in\left[0,1 \right] $, $\lambda=\lambda^{*}$ and $\lambda=1-\lambda^{*}$ will suffice. Due to the symmetry of $g$, branches are also created at the points $(1-\lambda^{*},0)$ and $(1-\lambda^{*},1)$ and the condition for these to meet is the same.
A remaining question is whether it is the middle eigenvalue that is equal to $1$. This is answered in a similar fashion to \cite{8} and hence, we shall require that
\begin{equation}
\textrm{tr}U^{2}_{A}-\det U^{2}_{A}-2+(\lambda^2-\lambda)\left|b_{AB}\right|^{2}+(\Lambda^2-\Lambda)\left|b_{AA'}\right|^{2}
\label{eq:41}
\end{equation}
\[+(\lambda^2-\lambda)(\Lambda^2-\Lambda)\eta^{2}\left|b_{AB}\right|^{2}\geq0\]
for all pairs $(\lambda,\Lambda)$ that make the interface possible. This will also imply that the other eigenvalues of $M^{T}M$, $\lambda_{1}$ and $\lambda_{3}$, are bounded away from $1$, i.e. $0\leq\lambda_{1}\leq\ 1 \leq\lambda_{3}$.

However, the formulae for the cubic-to-orthorhombic transformation get too involved and for this reason we will proceed numerically.
\subsection{Numerical Results} 
\begin{figure}[h]
	\centering
		\includegraphics[scale=1]{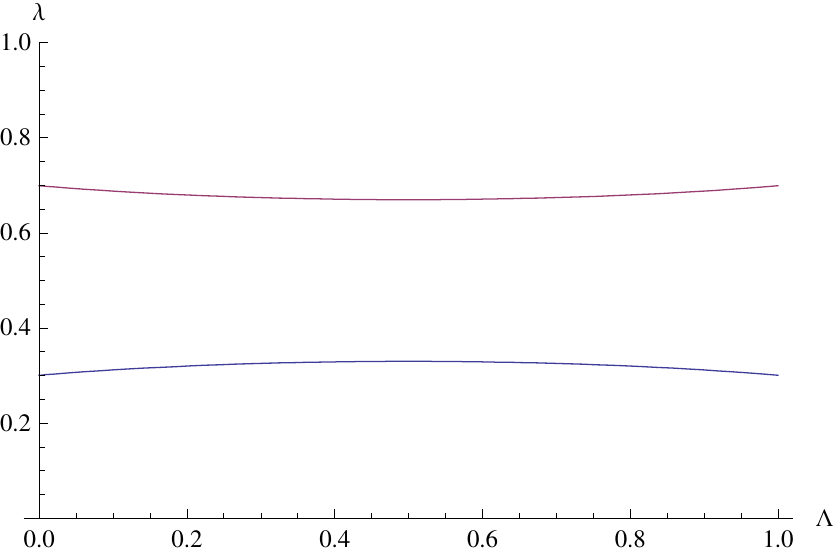}
	\caption{Values of $\lambda$ that make the interface possible for $0\leq\Lambda\leq1$}
	\label{fig:contourplot}
\end{figure}

In the remainder of this section we present a numerical calculation where, in accordance with \cite{9}, the lattice parameters for CuAlNi were chosen to be $\alpha=1.06372$, $\beta=0.91542$ and $\gamma=1.02368$. The martensitic variants used are the ones obtained from the experimental observations; that is, $A=3$, $B=6$ with $A'=4$, $B'=5$ their compound counterparts. Following the above analysis, we calculated the zeros of $g$, i.e. the values of $(\lambda,\Lambda)$ that allow the interface to occur. Relations (\ref{eq:36}) and (\ref{eq:38})-(\ref{eq:41}) were satisfied and we plotted $\lambda$ against $\Lambda$ as shown in Fig.~\ref{fig:contourplot}.

It is easily seen that the value of $\lambda$ does not change significantly from $\lambda^{*}$ (corresponding to $\Lambda=0$), which would give the classical interface between the Type-II twinning system and the austenite.

\begin{figure}[h]
\centering
\includegraphics[scale=0.7]{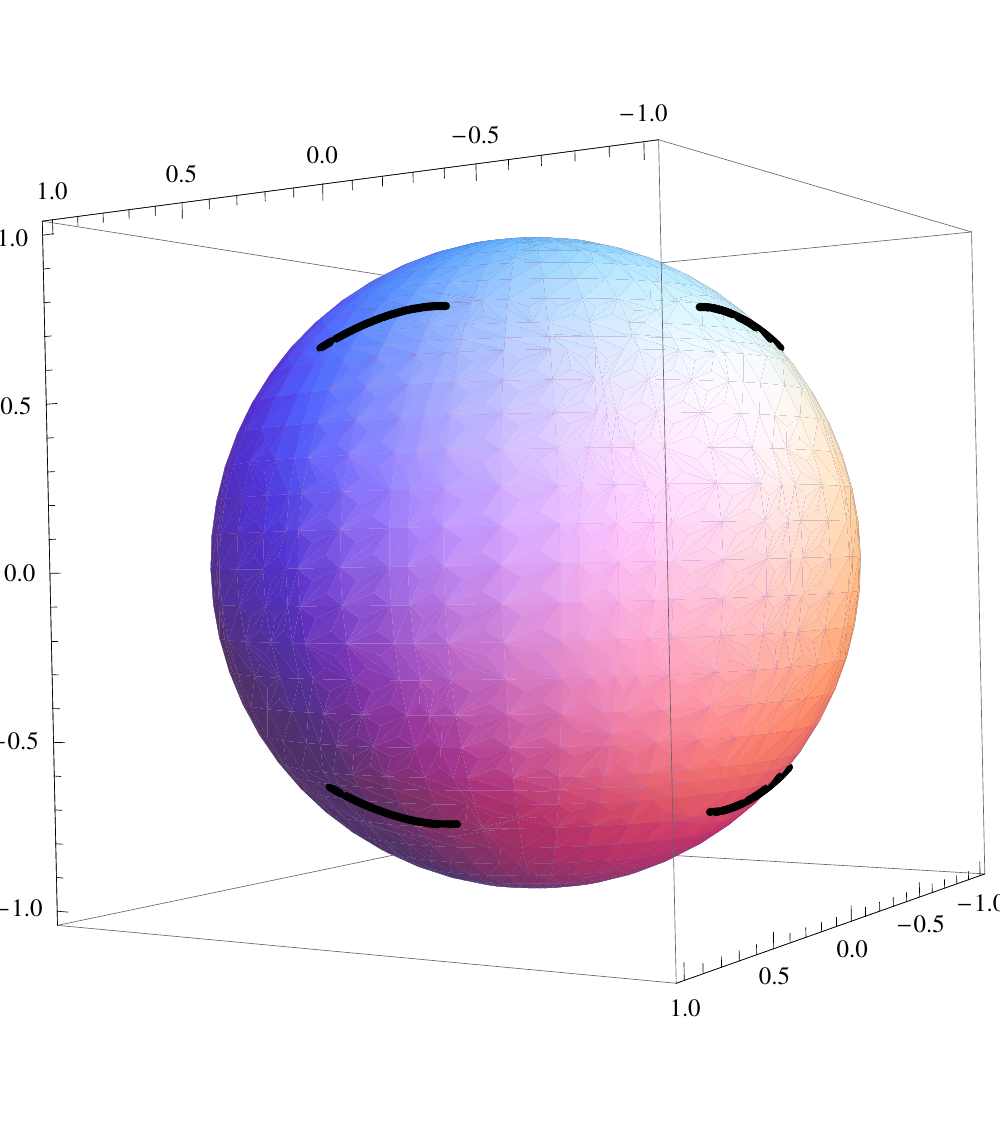} 
\caption{Plot of the normals on the unit sphere for $(\lambda,\Lambda)$ as in Fig.~\ref{fig:contourplot}}
\label{fig:NormalsOnSphere}
\end{figure}

Moreover, using the algebraic procedure given in \cite{8}, we calculated the different normals $m(\lambda,\Lambda)$ for $(\lambda,\Lambda)$ on the curves of Fig.~\ref{fig:contourplot} and a plot of these is given in Fig.~\ref{fig:NormalsOnSphere}, where the normals are depicted as points on the unit sphere.
These normals lie on four segments of curves whose endpoints correspond to the normals of possible classical austenite-martensite interfaces. A similar calculation was performed in \cite{1} for the cubic-to-tetragonal transformation, where, in contrast to the predictions here, these segments were in fact arcs of circles on the unit sphere.

A more detailed comparison of these theoretical predictions with the observed non-classical interfaces will appear in \cite{3}.

\section{Conclusion} This paper brings a theoretical analysis of compatible interfaces between austenite and a crossing twins microstructure of 2H martensite of the CuAlNi shape memory alloy, where the crossing twins microstructure consists of Type-II and compound twinning systems. These interfaces were recently observed by optical microscopy during the shape recovery process of single crystals of this alloy \cite{3}. The main aim of this paper is to show that these interfaces (although never observed in any shape memory alloy before) do not contradict the commonly accepted non-linear elasticity model, but are, on the contrary, predictable by this model for arbitrary volume fraction of the compound laminate. Since the relation (\ref{eq:37}) between this volume fraction $\Lambda$ and a volume fraction of the Type-II laminate $\lambda$ intersecting compatibly the compound twins was derived analytically for a general cubic-to-orthorhombic transition, the analysis brought by this paper can be easily applied to predict the existence of similar non-classical interfaces in any other shape memory alloy with the same class of transition (e.g. CuAlMn) as well as for materials undergoing the cubic-to-tetragonal transitions (since the tetragonal symmetry is a member of the orthorhombic symmetry class.) 

The numerical simulations carried out in the last subsection of this paper reveal that there is a dependence between the compound volume fraction and the habit plane orientation for the CuAlNi alloy. This finding is consistent with the optical observations of slightly curved non-classical interfaces between austenite and the crossing-twins microstructure with heterogeneous compound volume fraction (Fig.~\ref{Curved Interface}).

\section{Acknowledgements} The experimental part of this paper was financially supported by the project No. 202/09/P164 of the Czech Science Foundation and by the institutional project of IT ASCR v.v.i. (CEZ:AV0Z20760514). Both supports are acknowledged by H.S.. The theoretical part of the paper (J.M.B. and K.K.) was supported by the EPSRC New Frontiers in the Mathematics of Solids (OxMOS) programme (EP/D048400/1) and the EPSRC Science and Innovation award to the Oxford Centre for Nonlinear PDE (EP/E035027/1). The project originated in the EU MULTIMAT network (MRTN-CT-2004-505226).
\end{document}